%% file: paper.tex
\pdfoutput=1
\documentclass[a4paper,fleqn]{cas-dc}

\usepackage[numbers]{natbib}
\usepackage{bm}
\usepackage{siunitx}

\graphicspath{{./fig/}}
\include{macros}

\def\tsc#1{\csdef{#1}{\textsc{\lowercase{#1}}\xspace}}
\tsc{WGM}
\tsc{QE}
\tsc{EP}
\tsc{PMS}
\tsc{BEC}
\tsc{DE}

\begin{document}
\let\WriteBookmarks\relax
\renewcommand{\floatpagefraction}{0.7}
\renewcommand{\textfraction}{0.15}
\renewcommand{\topfraction}{0.85}
\renewcommand{\bottomfraction}{0.65}
\setcounter{totalnumber}{4}
\setcounter{topnumber}{3}
\setcounter{bottomnumber}{2}
\shorttitle{Adaptive Multiphysics Coupling}
\shortauthors{S. Candelaresi et~al.}

\title [mode = title]{Adaptive Multiphysics Coupling for Hyperbolic Systems}

\author[1,2]{S. Candelaresi}[orcid=0000-0002-7666-8504]
\ead{simon.candelaresi@uni-a.de, simon.candelaresi@gmail.com}

\author[3]{E. Faulhaber}[orcid=0000-0001-9788-5949]

\author[1,2]{M. Schlottke-Lakemper}[orcid=0000-0002-3195-2536]

\affiliation[1]{organization={High-Performance Scientific Computing, Institute of Mathematics, University of Augsburg},
                addressline={Universit\"atsstra\ss e 12a},
                postcode={86159},
                city={Augsburg},
                country={Germany}}

\affiliation[2]{organization={Centre for Advanced Analytics and Predictive Sciences, University of Augsburg},
                addressline={Universit\"atsstra\ss e 12a},
                postcode={86159},
                city={Augsburg},
                country={Germany}}

\affiliation[3]{organization={Department of Mathematics and Computer Science, University of Cologne},
                addressline={Weyertal 86-90},
                postcode={50931},
                city={Cologne},
                country={Germany}}

\begin{abstract}
For the discontinuous Galerkin code Trixi.jl we implement capabilities for adaptively
coupling arbitrarily many domains with different physics.
The coupling of the systems is realized through the exchange of boundary information
and user-definable coupling functions.
This gives us the ability to couple systems that do not share a single variable
or have a very different number of variables.
This is particularly useful when we have a hierarchy of systems.
Our implementation is such that we can adaptively select the model of the coupled systems
in the course of the simulation.
This allows for highly dynamical scenarios as can be found e.g.\ in astrophysics.
The criteria for adaptivity can be user defined and tailored to
the physical problem.
Compared to computing the most complex model in the entire simulation domain,
our method of coupled systems of high and low complexity
leads to a significant reduction in computational time with very little overhead.
\end{abstract}


\maketitle

\section{Introduction}

Complex real world problems can surprisingly often be described using
a mathematical model; why this is, and where it is limited,
has been subject to scientific and philosophic debate
(e.g.\ \cite{Wigner1960, Tegmark2008, hossenfelder2018lost}).
In computational fluid dynamics (CFD) we represent a fluid, gas or plasma
in terms of a set of partial differential equations defined
on a domain.
Since the systems are in practice dynamic we are dealing with the time
variation of the dependent variables.
Boundaries are in most cases also part of the model.
Mathematically we are then dealing with an initial-boundary value problem.

There exists a plethora of physical models in CFD,
as well as a multitude of numerical methods to solve these equations.
The result is a forest of numerical codes, each with its own
ideal use cases and problems for which it is unsuited.
The more complex numerical codes make use of different
numerical methods, depending on the equations that are being used
and the physical parameters.
This compartmentalization allows us to simulate more complex systems
and systems that employ sets of equations that are not practically solvable
using only one method.
The usage of different numerical methods and physical equations in
the same problem is sometimes referred to as {\em multiphysics}.

While single problems using one numerical method can be solved
relatively easily when separate, a {\em multiphysics} problem
requires some form of coupling between the systems.
How to design such a coupling is highly non-trivial,
which is why the preferred method is to avoid it altogether
and describe the coupled systems as one larger system that is being
solved using one physical model and one numerical method.

Multiphysics coupling is not an entirely new problem, as has been
pointed out by \cite{Dugeai-Girodroux-Lavigne-2-1-AeroLab-2011}
who give a short overview on coupling software.
Albeit cursory, it gives some insights into different physics that have been coupled.
The focus is on boundary coupling where the boundary values are simply exchanged,
which is similar to what we are aiming at.
An analogous approach is used by \cite{Francois-Xavier-311-CompSci-2008} for the heat equation.
They require that the heat values, together with the heat fluxes, on both sides of
the coupling interface are the same.
This is enforced through an iterative process.
With a heat value given in domain 1 they compute the flux from domain 1.
This can be used to calculate the value in domain 2 and with that the flux from domain 2.
The flux from domain 2 gives the value in domain 1.
This is repeated until it converges, if it converges.

On the technical side, there are aspects like parallelism that can
be challenging, particularly if we are trying to devise
an efficient code \citep{10.1007/978-3-319-46735-1_6}.
When coupling very different systems, we would like to avoid
idle computational resources.
This load balancing issue can have a significant effect on the simulations
\citep{10.1007/978-3-030-39181-2_15}, but is not addressed in this work.

Here we present an implementation of coupling of $n$ systems
(sometimes referred to as {\rm participants}) through an interface
for the numerical code
Trixi.jl \citep{schlottkelakemper2020trixi, schlottkelakemper2021purely, ranocha2022adaptive}.
We show how we can couple the systems through the exchange
of boundary information and bespoke coupling functions.
Furthermore, the coupled domains can be adaptive
if required.

\section{Problem Description}

The problems we are tackling here are those that can be described
in conservation form.
In the most general form we have $n$ systems (participants).
Each has its own set of unknowns,
i.e.\ dependent variables, $\uu_i$ that depend on space $\xx$ and time $t$.
The systems can be written in conservation form as
\begin{eqnarray}
\partial_t \uu_i + \nab\cdot \ff_i(\uu_i) = 0 \\
\uu_i(\xx, t = 0) = \uu_i^0 \\
\uu_i(\partial V_i, t) = \uu_i^{\rm BV},
\end{eqnarray}
where $\uu_i^0$ is the initial condition and $\uu_i(\partial V_i, t)$ the boundary values
on the respective domains $V_i$.
Applying a numerical solver on these independent systems, we would, in general,
use different methods, e.g.\ finite volume on one and finite difference on another.
Even the grid discretization does not need to be the same.

For simplicity we will limit most of our theoretical discussion to the coupling
of two systems $a$ and $b$ with variables $\uu_a$ and $\uu_b$.
These two systems can be coupled in different ways.
Here we will use coupling through the interface boundary values.
One example would be the heat coupling of a hot slab heated from below
and a fluid on top.
Another is a magnetohydrodynamic (MHD) system coupled with an Euler system.
Each of these initial boundary value problems can be tackled
independently using the appropriate numerical method.

In general, not all quantities that are present in system $a$ are present in system $b$
or vice versa.
For instance, if system $a$ was described using the MHD equations and system $b$ using the
Euler equations we would share the density, momentum and internal energy variables,
but not the magnetic field.
Therefore, any coupling through a boundary exchange has to make sure
the variables are correctly transformed.
More complex couplings are those that do not directly share any variables.
The Vlasov equations have as unknowns density distributions.
Those can be easily transformed into fields, like velocity and magnetic field.
However, the transformation from fields into distributions requires
certain assumptions and is non-trivial.

\section{Mesh Views}

Mesh objects or structures in numerical codes contain the information
of the meshes, like size, geometry and space discretization.
They are typically static in the domains they cover.
Trixi.jl is no different in that regard.
Here we require adaptive and dynamical domains that can change in time.
In order to achieve this and minimize the changes to the Trixi.jl code
that are not related to coupling, we introduce mesh views.
Those are Julia structures that work like ordinary mesh structures,
but are related to a parent mesh and point to the sub-mesh that belongs
to the corresponding sub-domain.
To the code they are treated the same as ordinary meshes throughout.

\begin{figure}\begin{center}
\includegraphics[width=\columnwidth]{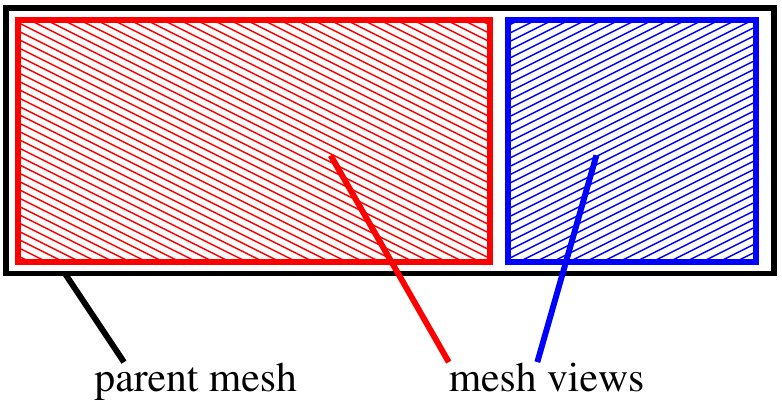}
\end{center}
\caption[]{
Parent mesh together with two mesh views acting like meshes on the subdomains.
}\label{fig: mesh_view}
\end{figure}

The related parent mesh holds the information
and data for the whole domain (see \Fig{fig: mesh_view}).
The mesh views are related to this parent mesh and hold information about their subdomain,
like subdomain size.
Although mesh views appear like meshes to the code, they are more flexible
and can be adapted during run time.

\section{Workflow}

The implementation of the coupling capabilities in Trixi.jl is such that
it is minimally invasive.
That is, any part that is not intended for coupling is left unchanged.
With that, any user who runs non-coupled simulations does not see
any difference and Trixi's infrastructure used in the simulation
thinks we are on a non-coupled system.
The mesh views are one example, as
they appear as standard meshes to Trixi.jl.

\begin{figure}\begin{center}
\includegraphics[width=\columnwidth]{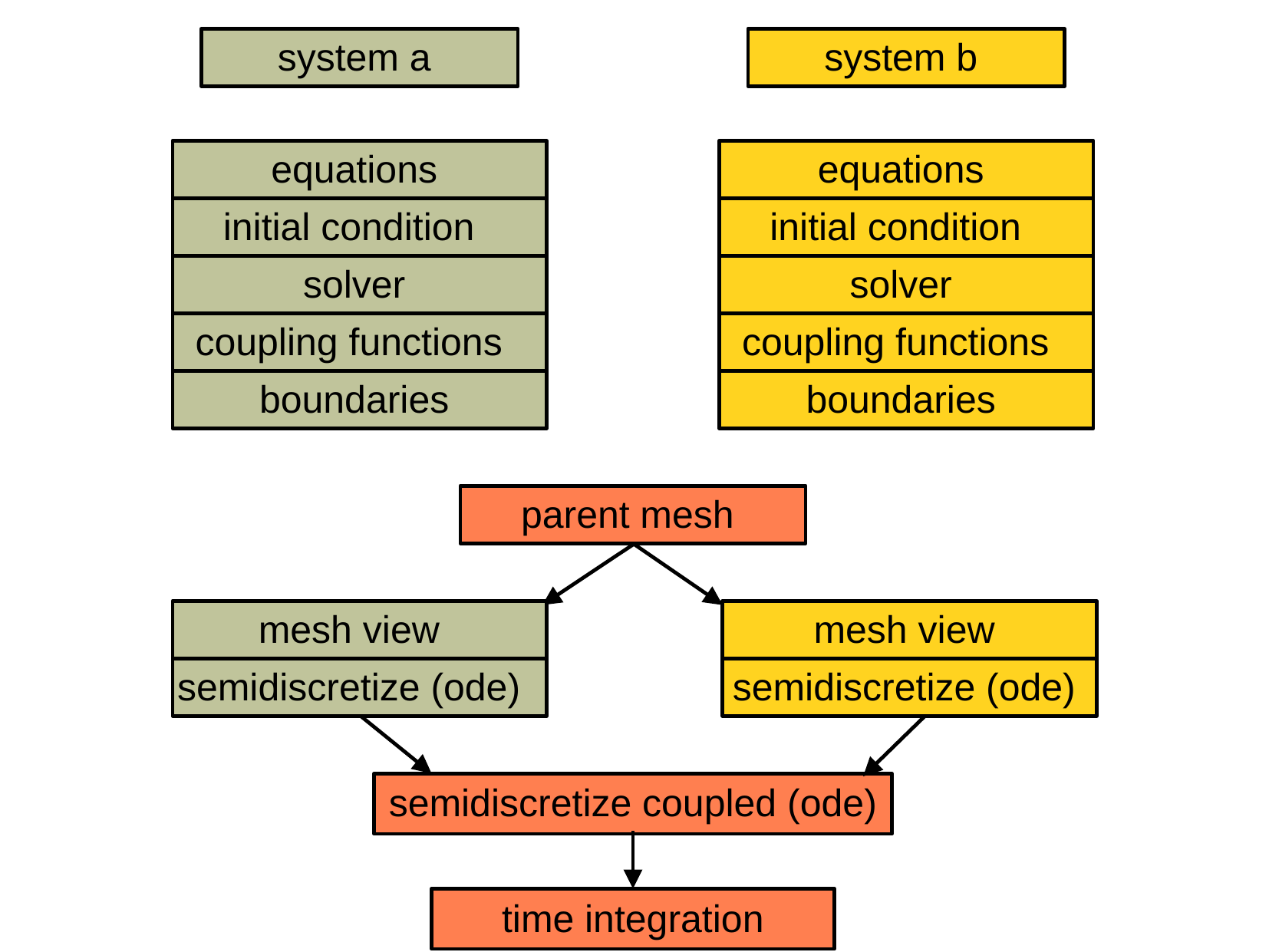}
\end{center}
\caption[]{
Outline of the workflow on how Trixi.jl handles coupled multiphysics simulations
for two systems.
}\label{fig: workflow}
\end{figure}

When setting up our coupled simulations we go through the same
steps as for a non-coupled one, but for each sub-system.
We define the problems with their initial conditions,
the solvers with their numerical fluxes and equations
(see \Fig{fig: workflow}).
Only when we define the boundary conditions we first need to define
coupling functions (see section \ref{sec: coupling functions}).

We then define the parent mesh on which we define as many mesh
views as we have sub-systems.
On each we apply Trixi's semidiscretization routine independently,
which transforms the PDE problem into an ODE problem.
Those semidiscretizations are then coupled in a coupled semidiscretization.
There, the coupling is being taken care of by exchanging
the coupled boundary conditions, while each system
is being solved separately.
Since the boundary values are required at every stage of the time integration,
our method couples more tightly than methods that couple every full time step.
A consequence is that the convergence order of the time integration
is preserved.

\section{Boundary Coupling}

Using the numerical code Trixi.jl \cite{schlottkelakemper2020trixi,
schlottkelakemper2021purely, ranocha2022adaptive} we
enhance the pre-existing boundary coupling
so that we are able to couple systems with more than
one variable, systems with different variables
and finally, make it adaptive.

All implementations are such that from a user's perspective there are
only a few changes compared to a non-coupled simulation.
Anyone not interested in coupling does not see any change in their code.
In its simplest form, this is done by implementing coupling boundary conditions that
copy the values from one domain to the boundary conditions
of the other domain.

\subsection{Coupling Functions}
\label{sec: coupling functions}

While a simple copy is sufficient for systems that share all variables,
like Euler and Navier-Stokes, we have to be more careful when coupling
systems with different variables.
For instance, MHD uses magnetic fields, while Euler does not.
Here we can simply copy the shared variables (density, momenta and energy)
and ignore the rest.

More complex coupling happens in hierarchical systems.
If we have a two-component gas in one domain and a one-component in the other
we could exchange the pressure, but not the density.
The two-component system has two densities.
A simple conversion would be to define two coupling functions,
one from system $a$ to $b$ ($c_{ba}$) and one from system $b$ to $a$ ($c_{ab}$):
\begin{eqnarray}
c_{ab} & = & c_{ab}(u_b) \\
c_{ba} & = & c_{ba}(u_a).
\end{eqnarray}
In general these functions depend on the equation parameters, such
as the adiabatic gas constant.
For a one-component gas (system $a$) coupled to a two-component gas (system $b$) in two dimensions we would have
the variables $\rho^{(a)}$, $v_1^{(a)}$, $v_2^{(a)}$, $p^{(a)}$ and $\rho_1^{(b)}$, $\rho_2^{(b)}$, $v_1^{(b)}$, $v_2^{(b)}$, $p^{(b)}$,
respectively.
Physically meaningful converter functions would be
\begin{eqnarray}
c_{ab} & = &
\begin{pmatrix}
\rho_1^{(b)} + \rho_2^{(b)} \\
v_1^{(b)} \\
v_2^{(b)} \\
p^{(b)}
\end{pmatrix}
\\
c_{ba} & = &
\begin{pmatrix}
\rho^{(a)}/2 \\
\rho^{(a)}/2 \\
v_1^{(a)} \\
v_2^{(a)} \\
p^{(a)}
\end{pmatrix}
.
\end{eqnarray}

\subsection{Coupling between Euler and Polytropic Equations}

A simple test example that is non-trivial is the coupling of the Euler equations
on one side and the polytropic equations on the other.
While the pressure is an independent variable for the Euler system,
it is related to the density in the polytropic system.
For the two-dimensional case we then have three variables in the polytropic
case and four in the Euler case.

The Euler equations in conservative form are
\begin{equation}\label{eq: Euler}
\partial_t
\begin{pmatrix}
\rho \\ \rho v_1 \\ \rho v_2 \\ \rho e
\end{pmatrix}
+
\partial_x
\begin{pmatrix}
 \rho v_1 \\ \rho v_1^2 + p \\ \rho v_1 v_2 \\ (\rho e +p) v_1
\end{pmatrix}
+
\partial_y
\begin{pmatrix}
\rho v_2 \\ \rho v_1 v_2 \\ \rho v_2^2 + p \\ (\rho e +p) v_2
\end{pmatrix}
=
0,
\end{equation}
with the density $\rho$, velocities $v_i$, energy density $e$ and pressure $p$.
Pressure and energy are related through
\begin{equation}\label{eq: p_rho_e}
p = (\gamma_\mathrm{E} - 1) \left(\rho e - \rho (v_1^2 + v_2^2)/2\right),
\end{equation}
with the adiabatic index $\gamma_\mathrm{E}$ which we choose to be $5/3$ for this test.

The polytropic equations are simpler and given as
\begin{equation}\label{eq: polytropic}
\partial_t
\begin{pmatrix}
\rho \\ \rho v_1 \\ \rho v_2
\end{pmatrix}
+
\partial_x
\begin{pmatrix}
 \rho v_1 \\ \rho v_1^2 + \rho^{\gamma_\mathrm{p}} \\ \rho v_1 v_2
\end{pmatrix}
+
\partial_y
\begin{pmatrix}
\rho v_2 \\ \rho v_1 v_2 \\ \rho v_2^2 + \rho^{\gamma_\mathrm{p}}
\end{pmatrix}
=
0
.
\end{equation}
Here we choose as adiabatic index $\gamma_\mathrm{p} = 2$.

To couple the Euler system to the polytropic system we can copy over the relevant variables:
\begin{equation}
c_{ba} = (\rho, \rho v_1, \rho v_2).
\end{equation}
The other way is trickier as we need to extract the pressure/energy from the density.
We use the adiabatic relation $p \propto \rho^{\gamma_\mathrm{E}}$ with the proportionality factor $1$,
so $p = \rho^{\gamma_\mathrm{E}}$.
Combined with the equation \eqref{eq: p_rho_e} we get
\begin{equation}
c_{ab} = \left(\rho, \rho v_1, \rho v_2, \frac{\rho^{\gamma_\mathrm{E}}}{\gamma_\mathrm{E}- 1} + \rho(v_1^2 + v_2^2)/2\right).
\end{equation}

For the domain we choose $[-2, 2]\times[-1, 1]$ with $64\times 32$ cells
and split it into two equal parts, with the polytropic system to the left
and Euler system to the right.
The $y$-boundaries are periodic, while the $x$-boundaries are coupled at $x = 0$
and $x = -2$ with $x = 2$.
That way we get a periodic topology.

\begin{figure}\begin{center}
\includegraphics[width=\columnwidth]{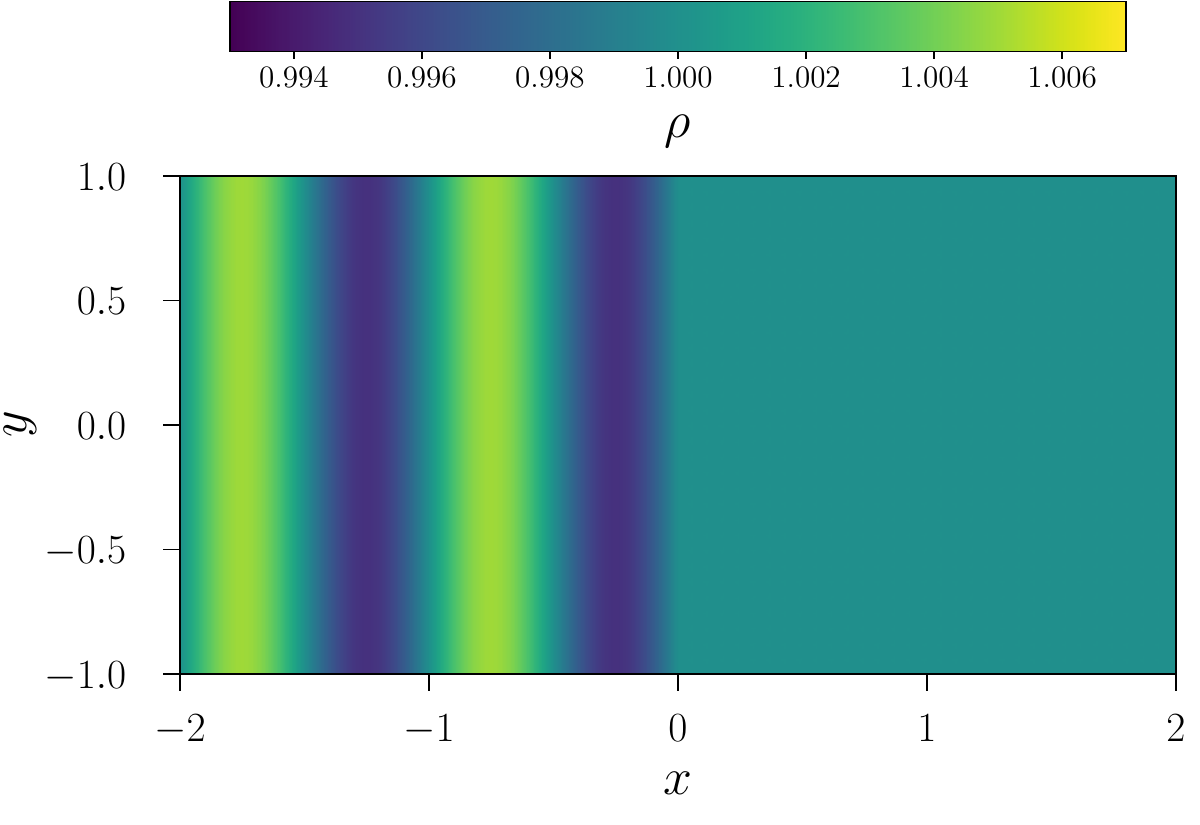}
\end{center}
\caption[]{
Density map for a polytropic system (left half) coupled to an Euler system (right half)
at time $0$.
}\label{fig: coupling_polytropic_euler_rho_t0}
\end{figure}

For the initial condition of the Euler domain we set the density to constant
$\rho = 1$, the velocities to $0$ and $p = 1$.
For the polytropic domain we prescribe a wave with positive traveling speed:
\begin{eqnarray}
\rho & = & (1.0 + 0.01 \sin(2\pi x))^{1/\gamma_\mathrm{E}} \\
v_1 & = & -0.01\sin(2\pi(x-1/2)) \\
v_2 & = & 0.0.
\end{eqnarray}
The initial density distribution is shown in \Fig{fig: coupling_polytropic_euler_rho_t0}.

\begin{figure}\begin{center}
\includegraphics[width=\columnwidth]{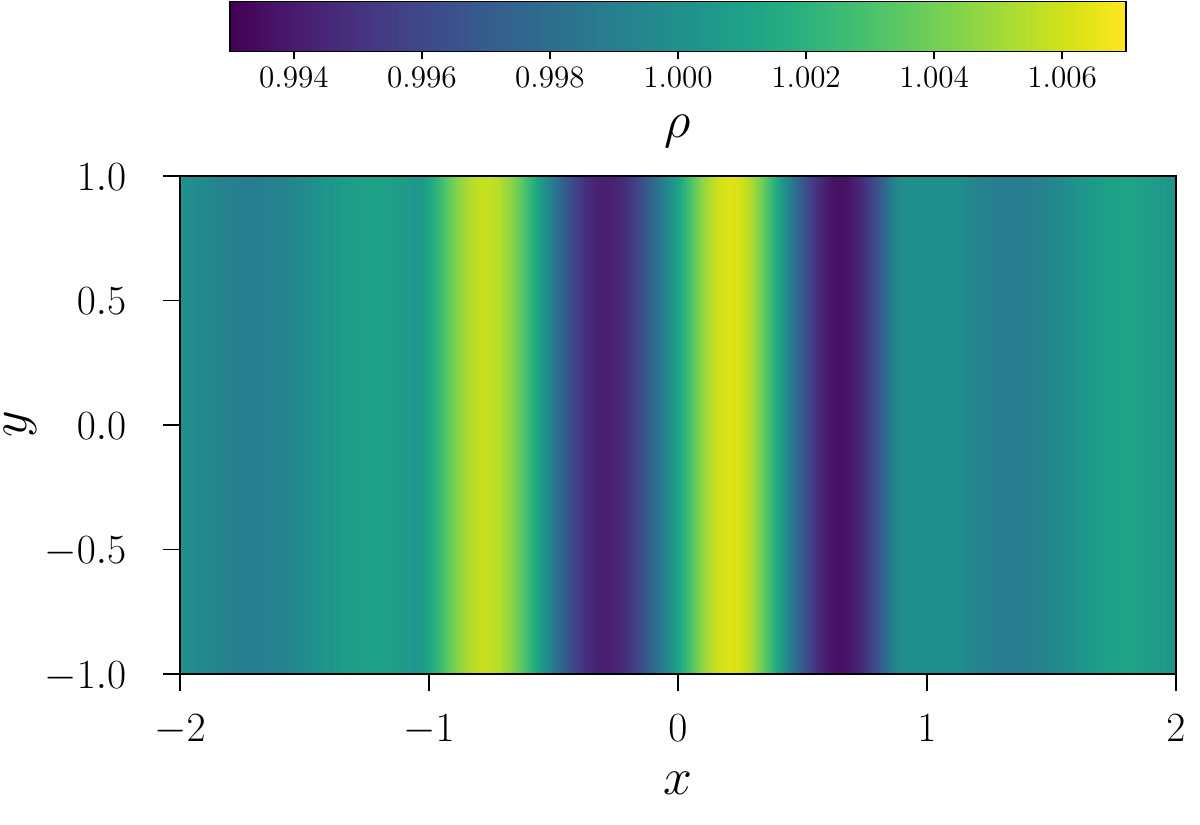}
\end{center}
\caption[]{
Density map for a polytropic system (left half) coupled to an Euler system (right half)
at time $0.6864$.
The different amplitudes at the two visible crests are a consequence of the different dispersion relations.
}\label{fig: coupling_polytropic_euler_rho_t0686}
\end{figure}

The wave travels from the polytropic domain into the Euler domain (\Fig{fig: coupling_polytropic_euler_rho_t0686}).
There, the dispersion relation changes and with that the shape of the wave.
We observe this in our simulations, as the wave's amplitude starts to change.
At the coupled boundaries we do not observe any discontinuities or any other
unphysical behavior.
This shows that our method and implementation are appropriate for the problem.

\subsection{Array of Coupled Systems}

\begin{figure}\begin{center}
\includegraphics[width=\columnwidth]{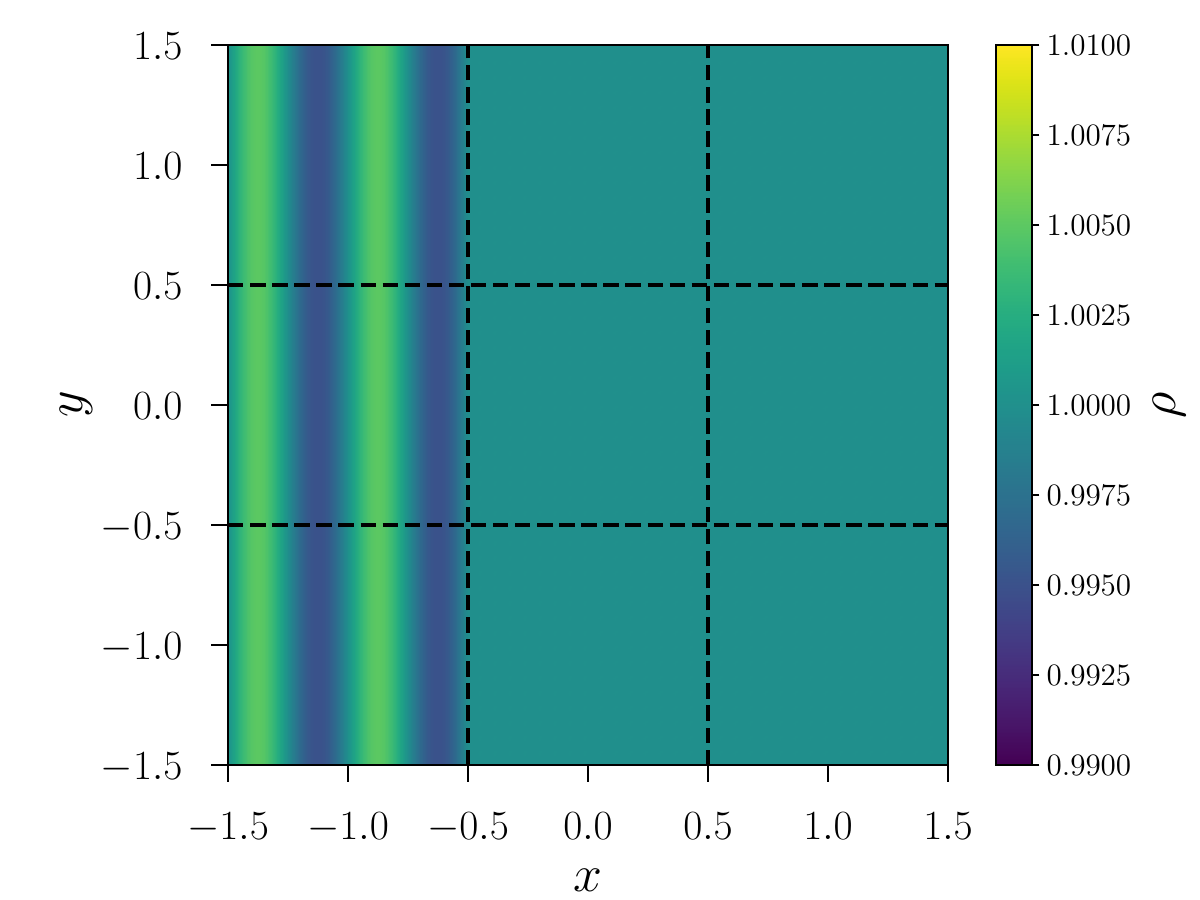}
\end{center}
\caption[]{
Initial density distribution for the array of coupled polytropic and isothermal systems.
The thick black dotted lines denote the inner coupling boundaries.
However, we also apply periodic coupling by coupling the top and bottom,
and the left and right.
}\label{fig: coupling_euler_euler_array_rho_t0}
\end{figure}

Due to the simplicity of the implementation we can easily couple multiple systems,
including entire arrays.
This is not limited to the $x$-coordinate, but can be just as easily done in the
$y$-direction.
To showcase this we perform simulations of $9$ systems arranged in a $3\times 3$ array.
For all systems we choose the polytropic equations \eqref{eq: polytropic}.
For the central system we choose $\gamma_\mathrm{p} = 1$ making it effectively isothermal.
The other $8$ systems have $\gamma_\mathrm{p} = 2$.
Initially there is a pressure wave given as
\begin{eqnarray}
\rho & = & \left\{ \begin{array}{l} (1 + 0.01\sin(4\pi x))^{1/\gamma_\mathrm{p}}, \ x < -0.5 \\ 1, \ x \ge -0.5 \end{array}\right.\\
\rho v_1 & = & \left\{ \begin{array}{l} (0.01\sin(4\pi(x-1/2))), \ x < -0.5 \\ 0, \ x \ge -0.5 \end{array}\right.\\
\rho v_2 & = & 0,
\end{eqnarray}
which corresponds to a right travelling pressure wave (\Fig{fig: coupling_euler_euler_array_rho_t0}).

For the coupling we use trivial coupling functions.
They are simply the identities $c_{ab} = u_b$ and $c_{ba} = u_a$.
To emulate a periodic domain we use this coupling through
the left and right, and top and bottom boundaries.

\begin{figure}\begin{center}
\includegraphics[width=\columnwidth]{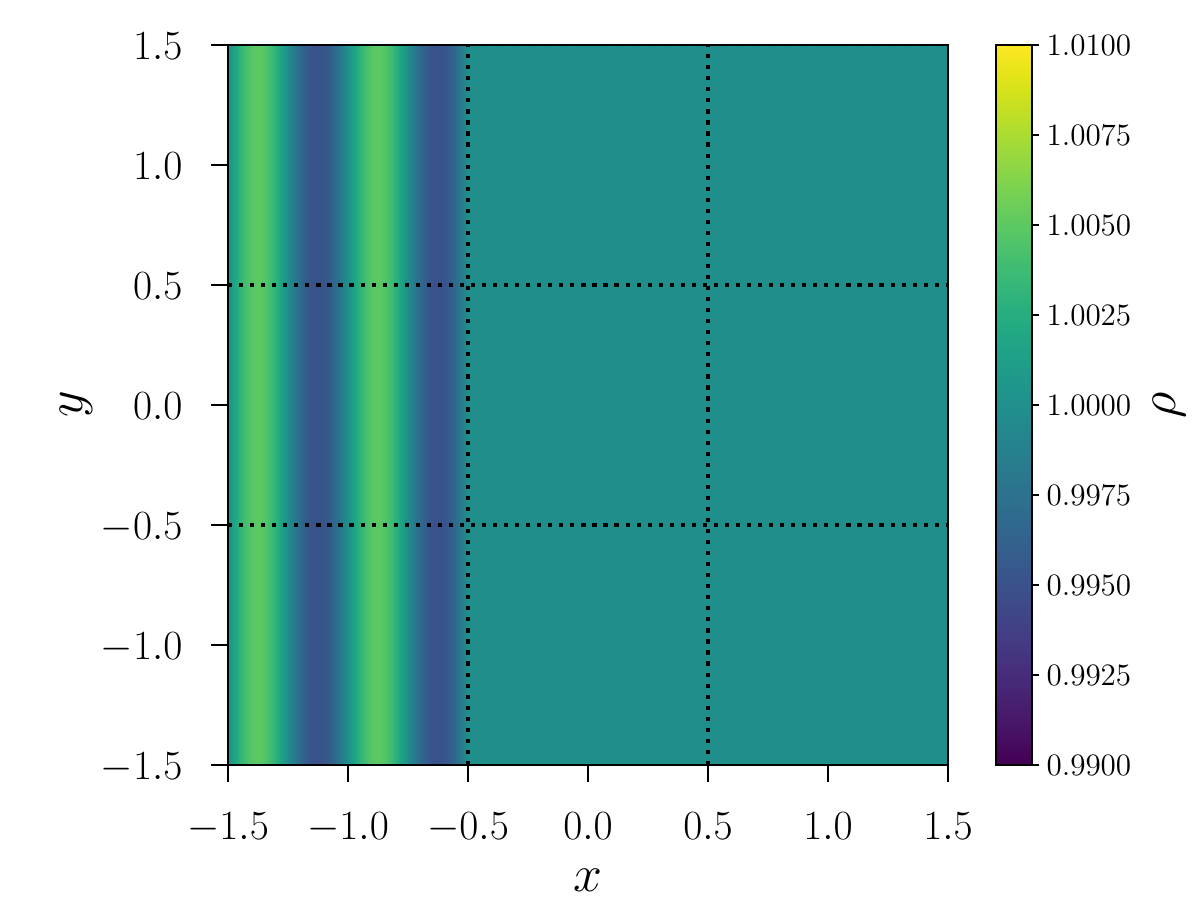}
\end{center}
\caption[]{
Density distribution at time $t = 0.823$ for a pressure wave traveling from polytropic
domains with $\gamma_\mathrm{p} = 2$ through an isothermal domain with $\gamma_\mathrm{p} = 1$ (center).
}\label{fig: coupling_euler_euler_array_rho_t0823}
\end{figure}

As the dispersion relation changes with $\gamma_\mathrm{p}$ we predict that the incoming
wave changes as it travels through the isothermal domain.
After time $t = 0.823$ we see a clear departure from a simple undisturbed linear wave
(\Fig{fig: coupling_euler_euler_array_rho_t0823}).
As with the coupled Euler and polytropic systems, we do not observe
any discontinuities at the coupled boundaries or any
other unphysical effects.

\section{Adaptive Model Selection (AMS)}

So far we kept our domains static and the models solved on each cell
are unchanged throughout the simulation.
This is useful for systems like a solid coupled to a gas.
But in fluid dynamics and astrophysics we are interested in scenarios
where the physics solved on the cells can change.
For instance, for a system with a strong magnetic field in
only parts of the domain we would choose to split the system so that
the MHD equations are solved there and the Navier-Stokes in the rest
of the physical domain to save some computational time.
But as the system is very dynamic, the magnetic field might travel
into the Navier-Stokes region, which would make it necessary
to adapt the coupling.

\begin{figure}\begin{center}
\includegraphics[width=0.8\columnwidth]{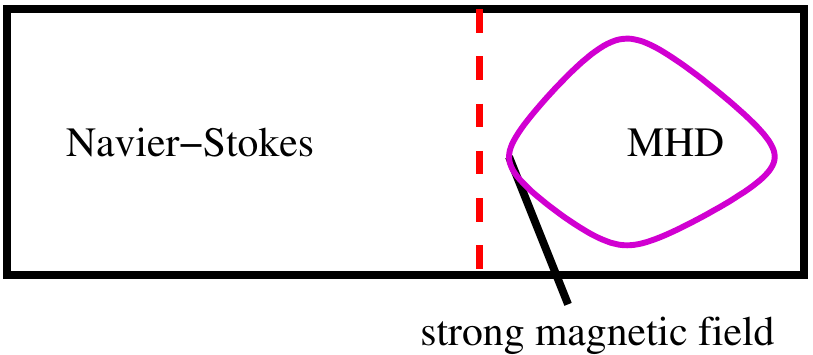} \\ \vspace{1em}
\includegraphics[width=0.8\columnwidth]{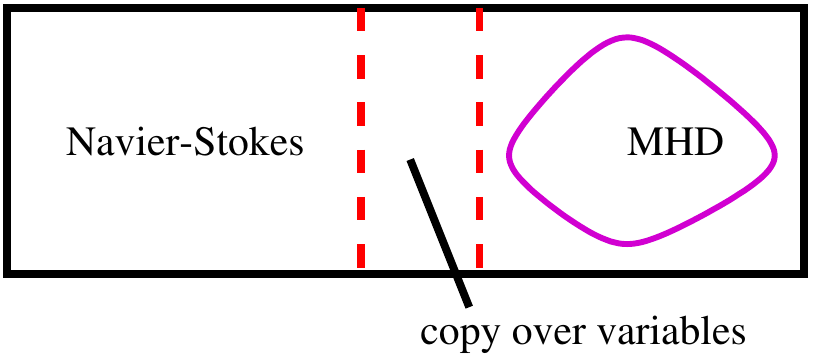} \\ \vspace{1em}
\includegraphics[width=0.8\columnwidth]{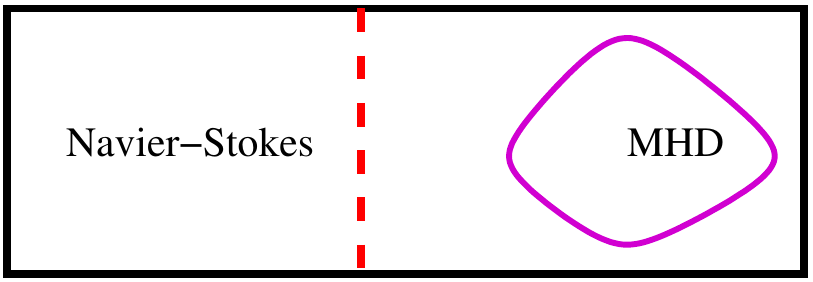}
\end{center}
\caption[]{
Sketch of the adaptive model selection process where we start with a Navier-Stokes system to
the left and an MHD system to the right.
In this scenario we assume that the magnetic field near the coupled boundary
becomes significant (top panel).
With our adaptive model selection, a layer within the Navier-Stokes domain is then
appropriated to simulate the MHD part, for which variables need to be copied (central panel).
Finally, we end up with an extended MHD domain (bottom panel).
}\label{fig: re_gridding}
\end{figure}

We do this by checking at every time step if we need to change
the models according to certain criteria.
Similar to adaptive mesh refinement (AMR), the user can define such criteria.
If the criterion for adaptive model selection is true
we change the mesh sizes of the sub-domains (mesh views).
This is much easier than if we were to change an actual mesh.
On the new resized meshes we then need to map the values of our
unknown vectors $\uu_i$ (\Fig{fig: re_gridding}).
Extending the MHD domain into the Navier-Stokes domain we would just take
the values of the density, velocities and pressure,
while setting the magnetic field to zero.
If we choose to extend the Navier-Stokes system we could
make the same copy of variables and ignore the magnetic field.
For the general case we could simply use the coupling functions,
which would make it consistent with the way we couple.

In general this works for any number of coupled boundaries.
But care has to be taken when defining the conditions.
When we extend the MHD system we need to set the magnetic field in
the reclaimed domain to $0$, as there was no magnetic field to copy
from Navier-Stokes.
If we then set the criterion for an extension of Navier-Stokes to be
how small the magnetic field is, we then would end up with
models switching back and forth at every time step.

After the values of the solution vector have been correctly copied
into the new solution vectors we construct a solution
vector for the combined system.
The size of this is in general different from the original
pre-AMS solution vector and the unknown variables are in general at
different indices.
Therefore, the problem that is being solved by the time integrator
changes and needs to be updated.

\subsection{Euler Equations Adaptively Coupled with Ideal MHD}

To showcase our adaptive model selection we
couple the Euler and ideal MHD equations using nine systems.
Eight of them are Euler and one is MHD.
The domain is $[-3, 3]\times[-3, 3]$.
The MHD domain is initially set at the center at the square $[-1, 1]\times[-1, 1]$,
while the Euler systems sit around it.
The coupling is done simply by exchanging the shared variables of density, momentum
and pressure while the magnetic field is of course not exchanged.

The MHD equations are the ones for a $2.5$-dimensional system.
This means that they are solved on a two-dimensional domain with
only $x$- and $y$-dependencies, but contain a $z$-component.
They are given in hyperbolic form as:
\begin{eqnarray}\label{eq: MHD}
& 0 =
\partial_t
\begin{pmatrix}
\rho \\ \rho v_1 \\ \rho v_2 \\ \rho v_3 \\ \rho e \\ B_1 \\ B_2 \\ B_3 \\ \Psi
\end{pmatrix}
+ \nonumber \\
 & \partial_x
\begin{pmatrix}
 \rho v_1 \\ \rho v_1^2 + p + E_\mathrm{mag} - B_1^2 \\ \rho v_1 v_2 - B_1 B_2 \\
 \rho v_1 v_3 - B_1 B_3 \\
 (E_\mathrm{kin} + \frac{\gamma}{\gamma - 1} p + 2E_\mathrm{mag})v_1 -
 B_1 v\cdot B + c_h\Psi B_1 \\
 c_h \Psi \\ v_1 B_2 - v_2 B_1 \\ v_1 B_3 - v_3 B_1 \\ c_h B_1
\end{pmatrix}
+ \nonumber \\
 & \partial_y
\begin{pmatrix}
 \rho v_2 \\ \rho v_1 v_2 - B_1 B_2 \\ \rho v_2^2 + p + E_\mathrm{mag} - B_2^2 \\
 \rho v_2 v_3 - B_2 B_3 \\
 (E_\mathrm{kin} + \frac{\gamma}{\gamma - 1} p + 2E_\mathrm{mag})v_2 -
 B_2 v\cdot B + c_h\Psi B_2 \\
 v_2 B_1 - v_1 B_2 \\ c_h \Psi \\ v_2 B_3 - v_3 B_2 \\ c_h B_2
\end{pmatrix},
\end{eqnarray}
with the magnetic field components $B_i$,
magnetic energy density $E_\mathrm{mag} = B^2/2$,
kinetic energy density $E_\mathrm{kin} = \rho v^2/2$, divergence cleaning field $\Psi$
and divergence cleaning speed $c_h$.
$\Psi$ is a Lagrange multiplier that, together with the last equation,
drives $\nabla\cdot B \to 0$.

\begin{figure}\begin{center}
\includegraphics[width=\columnwidth]{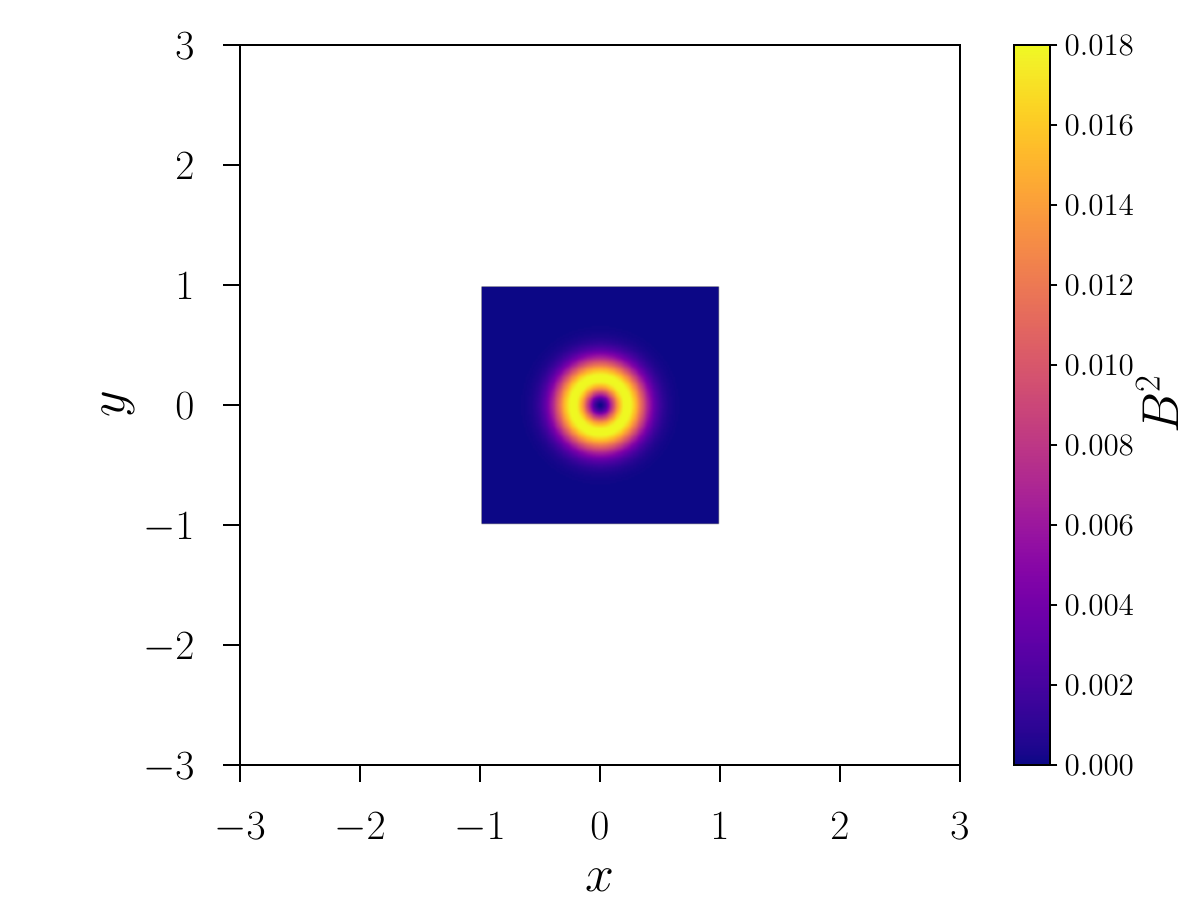}
\end{center}
\caption[]{
Magnetic energy density at initial time for the
adaptively coupled Euler-MHD system.
}\label{fig: coupling_euler_mhd_euler_adaptive_array_B2_t0}
\end{figure}

As initial condition we use constant density $\rho = 1$
and pressure $p = 1$.
We impose a constant velocity
by using the initial condition $\rho v_1 = 0.2$
and $\rho v_2 = 0.1$.
For the MHD part we choose a localized magnetic ring as initial condition.
This is expressed as
\begin{eqnarray}
B_1 = y e^{-10r^2}\qquad
B_2 = -x e^{-10r^2},
\end{eqnarray}
where $r = \sqrt{x^2 + y^2}$.
Since we are using a $2.5$-dimensional MHD model we have a $z$-component
that we choose to be $B_3 = 0$.
The generalized Lagrange multiplier we set to $\Psi = 0$ initially.
The initial magnetic energy density of this flux ring is shown in
\Fig{fig: coupling_euler_mhd_euler_adaptive_array_B2_t0}.

\begin{figure}\begin{center}
\includegraphics[width=\columnwidth]{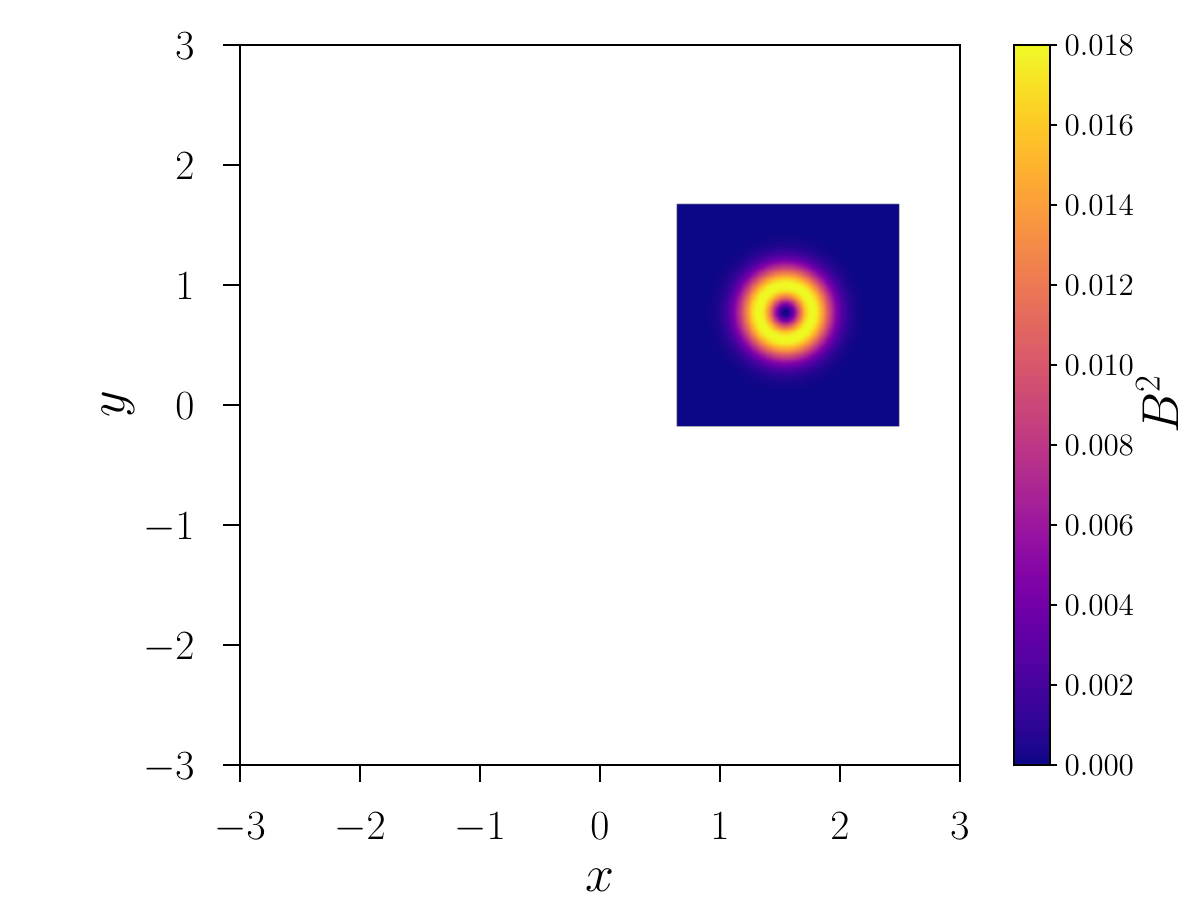}
\end{center}
\caption[]{
Magnetic energy density at time $t = 8$ for the
adaptively coupled Euler-MHD system.
}\label{fig: coupling_euler_mhd_euler_adaptive_array_B2_t8}
\end{figure}

In this test we use adaptive model selection where the domain for which we solve
the more computationally demanding MHD equations changes in time.
We define a criterion for when and how it changes.
If the magnetic field strength near the boundary of the domain
is above a threshold then the adjacent Euler element is converted
into an MHD element for the next time step.
If the magnetic field strength is lower than a different threshold, then
we convert the MHD cells into Euler cells.
As we populate new MHD cells with zero magnetic field we would run
into a constant switching back and forth situation.
To avoid this we impose a delay for which cells can be transformed
back to Euler cells after being transformed into MHD cells.
For this particular case an imposed delay of $t_\mathrm{delay} = 0.2$
gives the magnetic field enough time to propagate into the new cells.
The result is a model selection that smoothly moves with our magnetic ring
(\Fig{fig: coupling_euler_mhd_euler_adaptive_array_B2_t8}).

\begin{figure}\begin{center}
\includegraphics[width=\columnwidth]{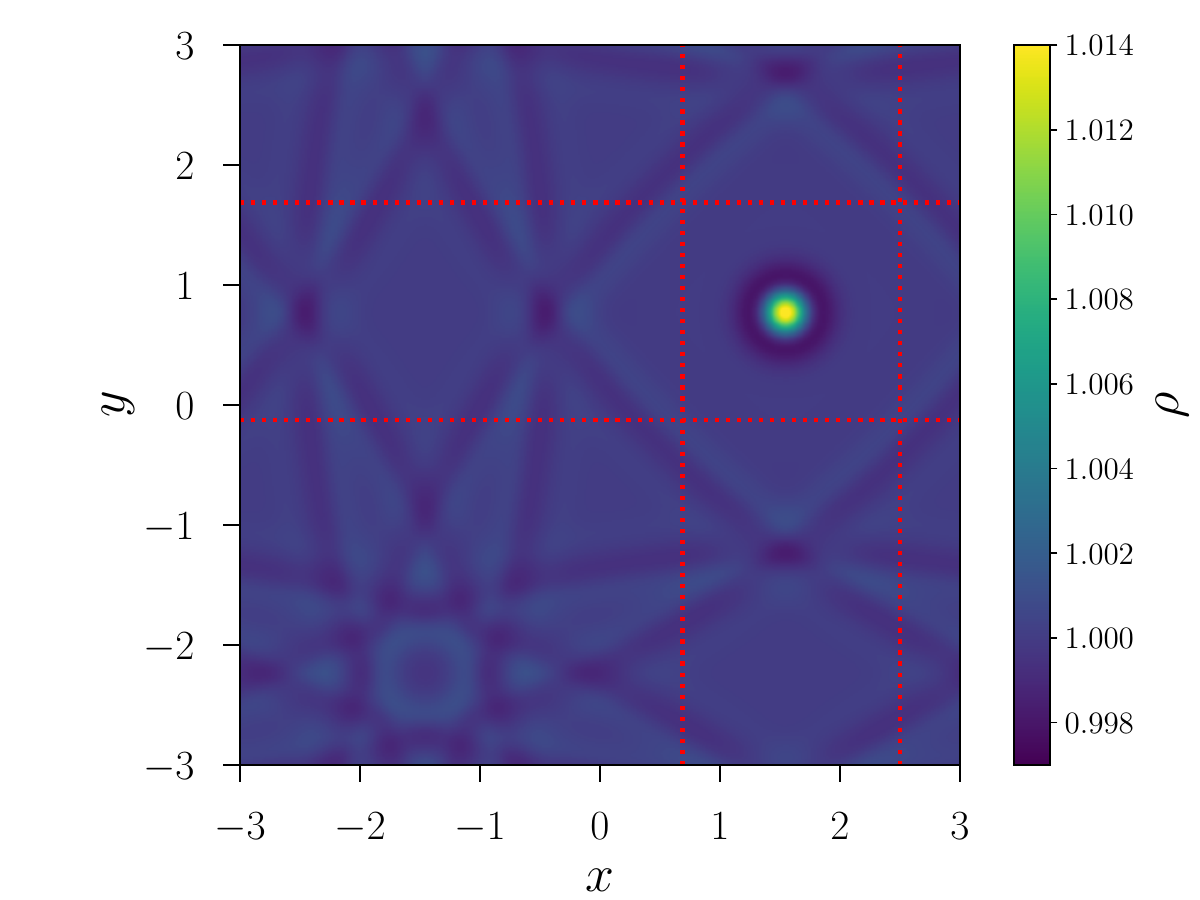}
\end{center}
\caption[]{
Density at time $t = 8$ for the adaptively coupled Euler-MHD system
with the domain boundaries as red dotted lines.
}\label{fig: coupling_euler_mhd_euler_adaptive_array_density_t8}
\end{figure}

From the magnetic field we cannot deduce if there are any
artificial effects, like discontinuities.
That is why we plot the density for all domains at time $t = 8$
(\Fig{fig: coupling_euler_mhd_euler_adaptive_array_density_t8}).
Since the field is initially not force-balanced,
the Lorentz forces accelerate the gas which leads to travelling waves.
During the entire evolution we do not notice any discontinuities
at the domain boundaries.

\begin{figure}\begin{center}
\includegraphics[width=\columnwidth]{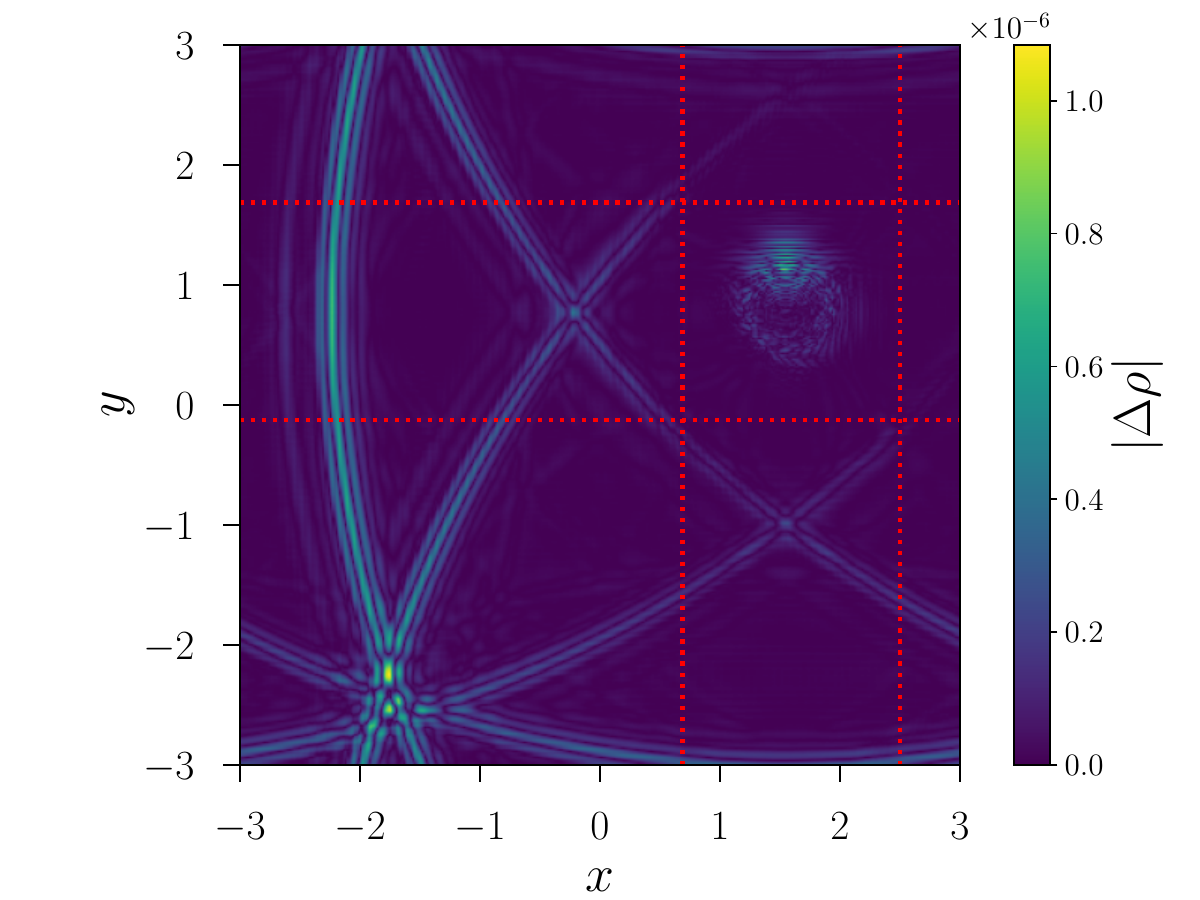}
\end{center}
\caption[]{
Unsigned difference of the density at time $t = 8$ between the adaptively coupled Euler-MHD system
and the system simulated using MHD throughout.
}\label{fig: coupling_euler_mhd_euler_adaptive_array_delta_density_t8}
\end{figure}

For comparison we also perform a simulation of the entire domain where we use the
MHD equations throughout.
This invariably leads to a somewhat different outcome.
If such differences are small for our application we can claim that
a coupled simulation produces similar results.
Here we compare the density of the system and find
that after time $t = 8$ the maximum difference is less than $1.1\times 10^{-6}$
(\Fig{fig: coupling_euler_mhd_euler_adaptive_array_delta_density_t8}).

\subsection{Computational Efficiency}

If it were not for performance, we would simulate the entire domain using
the more complex model.
In the above example we would choose to solve the MHD equations.
The solution of the Euler systems requires fewer calculations
for each element,
but we pay this reduction of computation with some overhead
in calculating the boundary conditions and changing the meshes.
Both require the evaluation of the coupling functions, which
in some cases can be complex.

Here we test the performance of the coupled Euler-MHD-Euler system and
compare it to a simulation where we solve the MHD equations
for the entire domain.
Since the system is two-dimensional and of low resolution we run
it without any parallelization on a laptop with an
Intel i7-11800H CPU operating at \SI{2.3}{\giga\hertz}.
Our Julia version is 1.10.6.
While the time steps $\dd t$ are of the same size, the computational time
per step is larger for the full-MHD system.
There, we require $\SI{133.5}{\ms}$ per step.
For the coupled simulation we require $\SI{44.9}{\ms}$.
This means there is almost a factor $3$ speed-up compared to
computing the complex model in the whole domain.
This is partially thanks to the relatively cheap copying
when computing the coupled boundaries.
The total time spent on that is only $0.8\%$ of the simulation time.
In the coupled simulation we spend ca.\ $\SI{1.9}{\ms}$ per time
step on checking if we should perform AMS.
This means $4.2\%$ of a time step is spent on this, which is larger
than the time we spent on copying the coupled boundaries.
To reduce the impact on the total wall time we could check only
every $n$ time steps instead.
While this happens at every time step, the actual AMS in our
adaptive simulations only happens when required.
This takes ca.\ $\SI{0.95}{\s}$ to perform.
Despite its overheads at every time step and especially when
performing AMS, our adaptively coupled simulation requires $\SI{756.4}{\s}$
to finish, while the MHD simulation takes $\SI{1896}{\s}$.
This corresponds to a $2.5\times$ speed-up in total wall-clock time, which is
smaller than the per-step factor of $3$ because of the AMS and coupling
overhead incurred at every time step.

\section{Conclusions}

We presented an implementation of adaptive multiphysics coupling for the
open source numerical solver Trixi.jl.
Its design makes it easy to use even for coupling of very different problems.
Using converter functions we can couple systems that do not share any variables,
which is of use for Vlasov systems coupled with MHD systems, among others.
By dynamically changing domains we can adapt our problem to changing situations,
like a traveling magnetic pulse.
The overhead from the coupling and domain adaptivity is of the order of
a few percent in the tested examples.
Thanks to the limited computational cost compared to using the most complex models
we can save a significant amount of computational time.
For the Euler-MHD system we save about two-thirds of the per-step computational cost,
corresponding to a $2.5\times$ reduction in total runtime,
compared to using MHD for the entire system.
The present work is a proof of concept on structured meshes -
extending the adaptive model selection to per-cell granularity and combining it with
adaptive mesh refinement on hierarchical Cartesian meshes is the natural next step.

\section*{Acknowledgements}

This project has benefited from funding by the Deutsche Forschungsgemeinschaft (DFG, German Research Foundation)
through the research unit FOR 5409 "Structure-Preserving Numerical Methods for Bulk- and Interface Coupling of
Heterogeneous Models (SNuBIC)" (project number 463312734).

\bibliographystyle{cas-model2-names}
\bibliography{references}

\end{document}

%% file: macros.tex
\newcommand{\EQ}{\begin{equation}}
\newcommand{\EN}{\end{equation}}
\newcommand{\EQA}{\begin{eqnarray}}
\newcommand{\ENA}{\end{eqnarray}}

\newcommand{\Fig}[1]{Figure~\ref{#1}}

{}
{}
{}

{}
{}
{}
{}
{}
{}
{}
{}
{}
{}
{}
{}
{}
{}
{}
{}
{}
{}
{}
{}
{}

{}
{}
{}

{}

{}
{}

\newcommand{\xx}{\bm{x}}

\newcommand{\uu}{\bm{u}}

\newcommand{\ff}{\bm{f}}

\newcommand{\nab}{\mbox{\boldmath $\nabla$} {}}

\newcommand{\dd}{{\rm d} {}}

\newcommand{\g}{\,{\rm g}}
\newcommand{\s}{\,{\rm s}}
\newcommand{\ms}{\,{\rm ms}}